\documentclass[a4paper,11pt]{amsart} 

\usepackage[a4paper,
            bindingoffset=0in,
            left=1in,
            right=1in,
            top=1.2in,
            bottom=1.1in,
            footskip=.25in]{geometry}

\usepackage{mathpazo}
\usepackage{euler}

\usepackage[normalem]{ulem}

\usepackage{setspace}
\usepackage[dvipsnames]{xcolor}

\usepackage{tikz}
\usepackage{tikz-cd}

\usepackage{mathtools}
\numberwithin{equation}{section}

\usepackage{graphicx}
\usepackage{enumerate} 
\usepackage{cancel}
\usepackage{comment}
\usepackage{latexsym}
\usepackage{amssymb,hyperref}
\hypersetup{colorlinks,linkcolor={MidnightBlue},citecolor={Maroon},urlcolor={OliveGreen}}
\usepackage[cp850]{inputenc}
\usepackage{epsfig}
\usepackage{amsthm}
\usepackage{amscd}
\usepackage{amsmath}
\usepackage{amsfonts}
\usepackage{graphics}
\newtheorem{theorem}{Theorem}[section]

\newtheorem{lemma}[theorem]{Lemma}
\newtheorem{problem}[theorem]{Problem}
\newtheorem{corollary}[theorem]{Corollary} 
 
\newtheorem{prop}[theorem]{Proposition}

\newtheorem*{theorem*}{Theorem} 
\newtheorem*{corollary*}{Corollary}
\theoremstyle{definition}
\newtheorem{example}[theorem]{Example}
\newtheorem{Counter-Example}[theorem]{Counter-Example}
\newtheorem{remark}[theorem]{Remark}
\newtheorem{definition}[theorem]{Definition}
\newtheorem*{remark*}{Remark}

\newtheorem*{definition*}{Definition}
\newtheorem*{example*}{Example}
\newtheorem{conj}[theorem]{Conjecture}

\newtheoremstyle{named}{}{}{\itshape}{}{\bfseries}{.}{.5em}{\thmnote{#3}}
\theoremstyle{named}

\newcommand{\BR}{\mathbb R} 
\newcommand{\BN}{\mathbb N} \newcommand{\BQ}{\mathbb Q}
 \newcommand{\BZ}{\mathbb Z}
\newcommand{\BF}{\mathbb F} 
 \newcommand{\BA}{\mathbb A}
\newcommand{\BP}{\mathbb P} 

\newcommand{\CE}{\mathcal E}

\newcommand{\CM}{\mathcal M} 
\newcommand{\CO}{\mathcal O}

\newcommand{\CU}{\mathcal U}

\newcommand\smvee{\raise0.3ex\hbox{$\scriptscriptstyle\vee$}}

\newcommand{\UX}{\underline{X}}
\newcommand{\UC}{\underline{C}}
\newcommand{\US}{\underline{S}}
\newcommand{\Uf}{\underline{f}}

\newcommand{\OCM}{\overline{\CM}}

\DeclareMathOperator{\Hom}{Hom}

\DeclareMathOperator{\rank}{rank}

\DeclareMathOperator{\Spec}{Spec}

\DeclareMathOperator{\bk}{\textbf{k}}
\DeclareMathOperator{\carac}{char}

\DeclareSymbolFont{greekletters}{OML}{cmr}{m}{it}
\DeclareMathSymbol{\vsigma}{\mathalpha}{greekletters}{"26}

\title{Campana Separable Rational Connectedness of Toric Orbifold}

\author{Enhao Feng}
\address{Department of Mathematics, University of Kansas, Lawrence, KS 66044, USA}
\email{fenge@ku.edu}
\author{Sara Mehidi}
\address{Mathematisches Institut
Heinrich Heine Universit\"at D\"usseldorf, Universit\"atsstr. 1
Geb\"aude 25.22.03.61, 40225 D\"usseldorf}
\email{sara.mehidi@hhu.de}

\begin{document}
\setstretch{1.1}

\begin{abstract}
    Campana introduced a notion of Campana rational connectedness for Campana orbifolds and conjectured that certain Fano orbifolds are Campana rationally connected. In this paper, we study a version of this conjecture in positive characteristic and verify a stronger statement for some toric orbifolds. We also provide examples illustrating the relation between the failure of separable Campana rational connectedness and singularities. 
\end{abstract}

\maketitle

\section{Introduction}

A \textit{rational curve} on a variety $X$ is the image of a non-constant morphism $f : \mathbb{P}^1 \to X$. The variety $X$ is said to be \textit{rationally connected} if roughly, any two general points can be joined by a rational curve. The geometry and arithmetic of a variety are largely influenced by the presence of such curves. However, in many natural geometric situations, one is led to consider pairs $(X,\Delta)$, where $\Delta$ is a boundary divisor; for instance, when compactifying a quasi-projective variety. In such a setting, the geometry depends partly on the behavior of curves relative to the boundary. This motivates the need for a refined notion of connectedness, in which one restricts the class of rational curves to the so-called \textit{Campana rational curves}, introduced by Campana  \cite[Section 6.1]{Cam04}.

\subsection{Campana orbifold}

Orbifold structures naturally arise in the study of rational points. For example, the distribution of rational points on a projective variety $Y$ can often be studied by considering a fibration $Y \to X$. Campana observed that endowing the base $X$ with additional data, namely, an orbifold structure (that is constructed from the fibration), provides a way to capture the arithmetic behavior of $Y$ (see \cite[Lecture 2]{Abr09}). In particular, he introduced the notion of an \textit{orbifold pair} in \cite{Cam04}, also called a \textit{Campana orbifold}.

A Campana orbifold is a pair $(X,\Delta_{\epsilon})$ of a projective variety $X$ and a $\BQ$-divisor
\[
    \qquad \Delta_{\epsilon} = \sum_i \left(1 - \frac{1}{m_i}\right) \Delta_i,
\]
where the $\Delta_i$ are distinct prime divisors, and $m_i \in \mathbb{N}_{>0} \cup \{\infty\}$ are multiplicities.
The divisor $\Delta_{\epsilon}$ records the presence of multiple fibres: given a fibration $f \colon Y \to X$, the base $X$ acquires a natural orbifold structure in which the multiplicities $m_i$ reflect the multiplicities of the fibres of $f$ along the divisors $\Delta_i$. In this way, the geometry of $Y$ is partially encoded in the pair $(X,\Delta_{\epsilon})$.

The orbifold canonical divisor is then defined by
\[
    K_{(X,\Delta_{\epsilon})} := K_X + \Delta_{\epsilon},
\]
and fundamental notions of birational geometry, such as Kodaira dimension, positivity, and the minimal model program, extend naturally to this setting.

\subsection{Campana rational connectedness}

The notion of rational connectedness admits a natural extension to the setting of Campana orbifolds (e.g. \cite{Cam04, Cam11}). Given a Campana orbifold $(X,\Delta_\epsilon)$, one considers Campana curves, i.e. curves whose tangency orders along the boundary divisors are compatible with the orbifold multiplicities. The resulting notion of \textit{Campana rational connectedness} requires roughly that any two general points of $X$ can be connected by a Campana curve. In positive characteristic, the appropriate notion is that of \textit{separable Campana rational connectedness}, where the family of rational curves making $X$ rationally connected retains nice deformation properties as in characteristic zero.

We are interested in the following problem:

\begin{problem}\label{prob: crc}
    Let $(X,\Delta_\epsilon)$ be a Campana orbifold over an algebraically closed field $\bk$ of arbitrary characteristic. When is $(X,\Delta_\epsilon)$ separably Campana rationally connected ?
\end{problem}

Campana studied a version of this problem and conjectured in \cite[Conjecture 6.27]{Cam11b} that a guiding principle of the theory is that the positivity of the orbifold anticanonical divisor $-K_{(X,\Delta_{\epsilon})}$ implies the existence of an abundance of Campana rational curves. The following formulation of the conjecture is from \cite[Conjecture 1.4]{CLT26}:
 
\begin{conj}\label{CampanaConj}
    Let $(X,\Delta_\epsilon)$ be a klt Fano orbifold over an algebraically closed field $\bk$ of characteristic zero; that is, $(X,\Delta_\epsilon)$ is a klt Campana orbifold ($m_i < \infty$ for all $i$) such that $-(K_X+\Delta_\epsilon)$ is ample. Then $(X,\Delta_\epsilon)$ is Campana rationally connected.
\end{conj}

On the other hand, one can distinguish two extreme cases of Problem \ref{prob: crc} depending on the coefficients of the boundary divisor $\Delta_\epsilon$. 

\begin{enumerate}[\hspace{5mm}1.]
    \item When $m_i = 1$ for all $i$, the boundary divisor $\Delta_\epsilon$ is trivial. We have the classical result that smooth Fano varieties are rationally connected in characteristic zero \cite{KMM92c, Cam94}, and Conjecture \ref{CampanaConj} expects this to be true in the orbifold setting.

    \item When $m_i = \infty$ for all $i$, the boundary divisor $\Delta_\epsilon$ is integral. In this case, the Campana orbifold $(X, \Delta_\epsilon)$ can be viewed as a compactification of the quasi-projective variety $X\backslash 
    \Delta_\epsilon$. Hence, the problem translates to the study of rational curves on quasi-projective varieties. 
\end{enumerate}

The second case will serve as a guiding example for us. In this paper, we allow some multiplicities to be infinity, so that $\Delta_{\epsilon}$ contains integral components. We refer to $(X, \Delta_\epsilon)$ in this case as a \textit{non-klt} Campana orbifold. In particular, we study Problem \ref{prob: crc} when the underlying variety $X$ is toric, which we refer to as a \textit{Campana toric orbifold}.

\subsection{Main results}

Our first result gives a sufficient condition for a non-klt smooth Campana toric orbifold to be separably Campana rationally connected.

\begin{theorem}[Theorem~\ref{General smooth case} and Corollary~\ref{3 adjacent smooth}]
\label{Intro General smooth case}
    Let $(X,\Delta_{\epsilon})$ be a Campana toric orbifold over an algebraically closed field $\bk$ of characteristic $p \geqslant \dim (X)$. Suppose there exists a smooth maximal cone in the fan defining $X$ such that all of its adjacent cones are smooth. Then $(X,\Delta_{\epsilon})$ is separably Campana rationally connected.
\end{theorem}

The klt case of this result (in all characteristics and without the assumption on the dimension) was previously established in \cite[Corollary 8.4]{CLT26}, whereas our approach extends to the non-klt setting. In particular, our dimension hypothesis is automatically satisfied when $X$ is a surface. Example \ref{dim=4<char=2} shows that, in higher dimension, the conclusion of the theorem may fail when this hypothesis is not satisfied. Although our result is stated in positive characteristic, it will be clear from the proofs that it also holds in characteristic zero, where separability is automatic. 

When $X$ is a toric surface, we further obtain a characterization of separability in terms of the singularities of the cones.

\begin{theorem}[Theorem~\ref{Crit_sing}]
\label{Intro Crit_sing}
    Let $(X,\Delta_{\epsilon})$ be a non-klt Campana toric surface orbifold over an algebraically closed field $\bk$ of characteristic $p>0$. If every ray of the fan of $X$ has a non-adjacent cone which is either smooth or tamely singular, then $(X,\Delta_{\epsilon})$ is separably Campana rationally connected.
\end{theorem}

Finally we also exhibit a family of Campana orbifolds for which the separability fails.

\begin{theorem}[Theorem~\ref{wps not sep}]
\label{Intro wps not sep}
    Let $\bk$ be an algebraically closed field of characteristic $p>0$. Let $(X,\Delta)$ be a non-klt Campana toric orbifold over $\bk$ such that $X$ is a weighted projective space $\mathbb{P}(q_0,\ldots,q_n)$ and every component of $\Delta$ has coefficient one. If $p$ divides one of the weights $q_j$, then $(X,\Delta)$ is not separably Campana rationally connected.
\end{theorem}

\subsection{The role of logarithmic geometry}
Logarithmic geometry allows one to keep track of intersection multiplicities along the boundary, so that tangency conditions are encoded directly in the structure of a morphism.

\begin{example}
    Let $f: \mathbb{P}^1 \to \mathbb{P}^2$ be a morphism such that $f^*\Delta_i = a_i P_i, \quad a_i \in \mathbb{Z}_{>0}$, where $\Delta_i = V(x_i)$ are the coordinate hyperplanes, and the points $P_1, P_2, P_3 \in \mathbb{P}^1$ are assumed to be distinct. We endow $\mathbb{P}^1$ with the logarithmic structure associated to the divisor $P_1 + P_2 + P_3$, and $\mathbb{P}^2$ with the logarithmic structure associated to the toric divisor $\Delta_1 + \Delta_2 + \Delta_3$. Then the morphism $f$ naturally upgrades to a morphism of log schemes. On the level of characteristic monoids of log structures, at each point $P_i$, the induced map is given by
    \[
        \overline{\mathcal{M}}_{\BP^2,f(P_i)} \cong \mathbb{N} \longrightarrow
        \overline{\mathcal{M}}_{\BP^1, P_i} \cong \mathbb{N},\quad 1 \longmapsto a_i.
    \]
\end{example}
Recently, the authors of \cite{CLT26} furnished Campana's orbifold theory within the framework of logarithmic geometry. 
Endowing $X$ with the logarithmic structure associated to $\Delta$ and $\mathbb{P}^1$ with an appropriate logarithmic structure encoding the marked points and contact orders, a morphism $\BP^1 \to X$ can naturally be lifted to a stable logarithmic map. In particular, families of Campana curves can be studied inside the well-behaved moduli space of stable log maps. This perspective provides a conceptual and flexible framework for studying the curves and their deformation theory \cite[\S 2,\S 3,\S 4]{CLT26}, which we have followed in this paper.

\subsection{Related work}
There are two notions that are closely related to Campana curves and Campana rational connectedness. The first one is \textit{orbifold rational curves} and \textit{orbifold rational connectedness}, which are developed by Campana in a series of work (\cite{Cam04}, \cite{Cam11}, etc) studying special varieties. This is the precursor to the theory adopted in this paper and it is shown in \cite{CLT26} that orbifold rational connectedness is equivalent to Campana rational connectedness. The second notion concerns $\BA^1$-curves and $\BA^1$-connectedness. This can be viewed as an analogue of rational curves on quasi-projective varieties, as suggested in \cite{KM99}, where the authors studied free rational curves on log Fano surfaces. Moreover, logarithmic geometry is also used as a framework to study $\BA^1$-connectedness (\cite{CZ14}, \cite{CZ16}, \cite{CZ19}, etc). For example, separable $\BA^1$-connectedness for log smooth varieties is studied in \cite{CZ19}, while \cite{CC21} explored its failure with respect to the strangeness of $\BA^1$-curves.

\vspace{3mm}

\noindent\textbf{Acknowledgments.} The authors are grateful to Qile Chen for proposing the question and many useful conversations regarding the problem. We thank Ana Maria Botero for insightful discussion on convex geometry and Marta Pieropan for feedback on a previous version. E. Feng would like to thank Brian Lehmann for helpful suggestions and all the encouragement and support he has given. The authors thank the referee for useful comments. S. Mehidi is supported by the NWO grant VI.Vidi.213.019.

\section{Background\label{Background}}
 
Throughout the paper, we work over an algebraically closed field $\bk$. Given a log scheme $X$ with log structure $\CM_X$, we denote by $\underline{X}$ its underlying scheme. Since we only use divisorial log structures in our results (\S \ref{Smooth Campana non-klt toric orbifold}, \S \ref{Singular Campana non-klt toric orbifold}), we recall their definition. We refer to \cite{katooriginal, ogusloggeom} for further details. If
\begin{itemize}
    \item $X$ is a regular scheme and $D$ a normal crossing divisor on $X$, or
    \item $X$ a toric variety and $D$ it's toric boundary
\end{itemize}
we let $j: U:= X-D \hookrightarrow X$. Then $\mathcal{M}_X:=\mathcal{O}_X \cap j_*O_U^{\times} \to \mathcal{O}_X$ defines a logarithmic structure on $X$. For example, given a curve $C/\bk$, a set of marked points induces a divisorial log structure on the curve. 

\subsection{Stable log map}
In this section, we briefly recall the definitions of stable log maps and contact order, which play a central role in the notion of Campana rational connectedness. Although these concepts arise from log geometry, in the toric setting they admit explicit combinatorial descriptions, which we use throughout this paper. We refer the interested reader to \cite[Sections 2 and 3]{CLT26} for additional details.

\begin{definition}
    A \textit{log curve} over a log scheme $S$ is a pair
    \[
        (\pi: C\to S, \{p_1,\dots, p_n\})
    \]
    satisfying the following conditions
    \begin{enumerate}[\hspace{5mm}1.]
        \item The underlying pair $(\underline{\pi}: \underline{C} \to \underline{S}, \{p_1,\dots,p_n\})$ is a family of pre-stable curves over $\US$ with $n$ marked points.
        \item $\pi$ is proper, log smooth, and integral.
        \item If $U\subset \UC$ is the smooth locus of $\underline{\pi}$, then $\CM_C|_U \cong \underline{\pi}^*\OCM_S \oplus \bigoplus_{k=1}^n p_{k_*}\BN_{\US}$.
    \end{enumerate}
\end{definition}

\begin{definition}
    A \textit{log map} over $X$ is a morphism of log schemes $f: C\to X$ such that $\pi:C\to S$ is a log curve over $S$. A log map is \textit{stable} if the underlying map $\Uf$ is stable in the usual sense.
\end{definition}

Given a stable log map $f:C\to X$, we say $f$ is \textit{non-degenerate} if $S$ is a log point with the trivial log structure. By \cite{katoF}, this implies that the underlying curve $\UC$ is an irreducible smooth curve and $\CM_C$ is the divisorial structure given by the marked points. Furthermore, given that the dense open subscheme of $\underline{C}$ where the log structure is trivial is sent by $\underline{f}$ to the dense open $X\backslash \Delta$ , we can conclude that $f(C)\not\subset \Delta$. In particular, $f^{-1}(\Delta) \subset \{p_1,...,p_n\}$.
 
\subsection{Contact order} 
Let $\underline{X}$ be a projective toric variety. Let $M$ be the lattice of characters and $N$ its dual lattice. We endow $\underline{X}$ with the divisorial log structure induced by its toric boundary $\Delta$, i.e. $\CM_X$ is the subsheaf of $\mathcal{O}_{\underline{X}}$ consisting of regular functions which are invertible on $\underline{X} \backslash \Delta$. We write $X$ for the induced log scheme. We let $\Sigma^{(1)}\subset \Sigma$ denote the collection of the rays of the fan $\Sigma$ of $\UX$. For each $\rho_i \in \Sigma^{(1)}$, let $u_{\rho_i} \in \rho_i \cap N$ denote its lattice generator and let $\Delta_i$ be the corresponding torus invariant divisor (we often abuse $\rho_i$ for its ray generator). We write $\Delta = \sum_{\rho_i\in \Sigma^{(1)}} \Delta_i$.

Let $f:C \to X$ be a non-degenerate stable log map with marked points $p_1,\dots, p_n$. The contact order at $p_i$ is given by
\[
    \overline{\CM}_{X,f(p_i)}\longrightarrow \overline{\CM}_{C,p_i} \cong \mathbb{N}.
\]
The contact order can be read from the fan of the toric variety and we will explain that. The point $f(p_i)$ lies in a unique torus orbit corresponding to some cone $\tau \in \Sigma$. The log structure of $X$ in a neighborhood of $f(p_i)$ is induced by the affine toric variety $X_\tau = \Spec \bk[M \cap \tau^{\vee}]$. Set $P_\tau = M \cap \tau^{\vee}$. We denote by $P_\tau^{\times}$ the subgroup of invertible elements of the monoid $P_\tau$. We have
\[
    \overline{\CM}_{X_{},f(p_i)}  =\overline{\CM}_{X_{\tau},f(p_i)} \cong P_\tau/P_\tau^{\times}.
\]
We deduce that:
\begin{align*}
    \Hom(\overline{\CM}_{X_{\tau},f(p_i)},\BN) & \cong \Hom(P_\tau/P_\tau^{\times},\BN)\\
    & \cong \Hom (P_\tau, \BN)\\
    & \cong N \cap \tau \subset N.
\end{align*}
More generally, the set of possible contact orders at a marked point is in bijection with the lattice points in the support of the fan of $X$, that is, the points on the cones which lie on the ambient lattice of the toric variety.

Contact orders are subject to intersection-theoretic constraints:

\begin{prop}[Balancing condition] \cite[Proposition 6.4.1]{CLS11}
   The collection of contact orders $\{c_{k}\}_{k=1}^n$ determined by $f:C\to X$ satisfies the balancing condition:
    \begin{equation}\label{balancing}
        \sum_{k=1}^n c_k = \sum_{i,k} c_{k,i}u_{\rho_i} = 0
    \end{equation} 
\end{prop}

\begin{example}
    Let $X$ be $\BP^m$ endowed with the log structure induced by its toric boundary $\Delta$ given by ${\rho_0},\dots, \rho_m \in \Sigma^{(1)}$. Let $H_i$ be the hyperplane section corresponding to $\rho_i$. We have $\sum_{i=0}^m u_{\rho_i} = 0$. Let $f:C\to \BP^m$ be a non-degenerate stable log map and let $\beta$ denote its curve class in $X$ and a collection of non-zero lattice points $\vsigma = \{c_k\in N\}_k$, where each $c_k$ specifies the contact order at the $k$-th marking $p_k$. By smoothness of $X$, if the coefficients outside the smallest cone containing $c_k$ are set to zero, there is a unique presentation 
    \[
        c_k = \sum_{i=0}^m c_{k,i}u_{\rho_i}.
    \]
    The balancing condition (\ref{balancing}) is a consequence of the intersection-theoretic constraints:
    \[
        \sum_{k=1}^{|\vsigma|} c_{k,i} = H_i\cdot \beta.
    \]
\end{example}

\subsection{Campana rational connectedness for toric orbifolds.}

In this section, we recall the notion of Campana rational connectedness. We first adapt the definition of a Campana orbifold (\cite[Definition 5.1]{CLT26}) to the toric case. 

\begin{definition}[Campana toric orbifold]\label{orbifold}
    Let $\UX$ be a projective toric variety with toric boundary $\Delta = \bigcup_i \Delta_i$. Let $X$ be the induced log scheme. To each irreducible component $\Delta_i$, we assign a weight 
    \[
        \epsilon_i = 1 - \frac{1}{m_i}
    \]
    where $m_i\in \BZ_{\geqslant 1}\cup\{\infty\}$. Define the $\BQ$-divisor
    \[
        \Delta_{\epsilon} = \sum_{i} \epsilon_i\Delta_i
    \]
    The pair $(X,\Delta_\epsilon)$ is called a \textit{Campana toric orbifold}. If $\epsilon_i\neq 1$ for all $i$, $(X,\Delta_\epsilon)$ is called a \textit{klt Campana toric orbifold}. Otherwise, it is called a non-klt Campana toric orbifold. If all the $m_i$'s are $\infty$, we say that it is an \textit{absolute non-klt Campana orbifold}. 
\end{definition}

\begin{definition}\cite[Definition 5.2]{CLT26} \label{Campana type} 
    Let $(X,\Delta_\epsilon)$ be a Campana toric orbifold. Let $f:C\to X$ be a stable log map over a log point $S$. Denote $c_{k,i}$ the multiplicity of the $i$-th component $\Delta_i$ along the $k$-th marked point, and let $\vsigma=\{c_{k,i}\}$ be the set of contact orders. We say that $\vsigma$ is of \textit{Campana type} if 
    \begin{enumerate}[\hspace{5mm}1.]
        \item When $m_i < \infty$, every index $k$ satisfies either $c_{k,i} = 0$ or $c_{k,i}\geqslant m_i$.
        \item When $m_i = \infty$, there is at most one index $k$ with $c_{k,i} > 0$.
    \end{enumerate}
\end{definition}

\begin{definition}[Campana rational connectedness]\cite[Definitions 4.1 and 5.6]{CLT26}\label{CRC}
    Let $X$ be the log scheme over $\bk$ induced by a projective toric variety of dimension $d$ with toric boundary $\Delta$. We say that $X$ is (separably) rationally $\vsigma$-connected if there is a family of genus zero non-degenerate stable log maps $\pi:\CU\to W$, $ev:\CU\to X$ with contact orders $\vsigma$ and such that $\underline{\CU} = W\times\BP^1$ and
    \[
        ev^{(2)}:\CU\times_W\CU = W\times \BP^1\times\BP^1 \to X\times X
    \]
    is dominant (and separable). Moreover, we say that $(X,\Delta_\epsilon)$ is \textit{(separably) Campana rationally connected} if there exists a collection of contact orders $\vsigma$ of Campana type such that $X$ is (separably) $\vsigma$-rationally connected. When all non-zero contact orders of $\vsigma$ are not divisible by $\carac \textbf{k}$, we say $(X,\Delta_\epsilon)$ is (separably) Campana rationally connected by good contact orders. 
\end{definition}

When $\bk$ is uncountable, rational connectedness is equivalent to the geometric condition that for a general pair of points $x, y \in X$, there is a log map $f:\BP^1\to X$ with contact order $\vsigma$ satisfying $x,y\in f(\BP^1)$. When $(X,\Delta_{\epsilon})$ is a Campana toric orbifold, the authors gave a sufficient combinatorial condition for Campana rational connectedness in terms of contact orders and the lattice.

\begin{theorem}
    \cite[c.f. Theorem 8.2]{CLT26} \label{QileMain}\label{CLT8.2} Let $\UX$ be a projective toric variety of dimension $d$ with toric boundary $\Delta$. Let $\vsigma$ be a collection of positive contact orders satisfying the balancing condition. Then  
    \begin{enumerate}
        \item If the sub-lattice $\BZ\cdot \vsigma\subset N$ is of rank $d$, $X$ is rationally $\vsigma$-connected.
        \item If $\BZ\cdot\vsigma$ is of rank $d$, and $N/(\BZ\cdot\vsigma)$ contains no $\carac \textbf{k}$-torsion, then $X$ is separably rationally $\vsigma$-connected.
    \end{enumerate}
\end{theorem}

Note that the original theorem is stated with $\UX$ being smooth and $\Delta$ being SNC, but as mentioned in \cite[Remark 8.3]{CLT26}, these assumptions are not necessary. In particular, any klt Campana toric orbifold is Campana rationally connected by good contact orders:

\begin{corollary}\cite[Corollary 8.4]{CLT26}\label{campana klt} 
    Endow the toric variety $\UX$ of Theorem \ref{CLT8.2} with the log structure induced by $\Delta$ and consider a Campana klt toric orbifold $(X,\Delta_{\epsilon})$. Then $(X,\Delta_{\epsilon})$ is separably Campana rationally connected by good contact orders.
\end{corollary}

\begin{remark}\label{remark on non klt CRC}
    As noted in \cite[Remark 5.3]{CLT26}, the definition of Campana rational connectedness in the non-klt case (Definition \ref{Campana type}) is weaker than the one present in the literature. When the multiplicities of a boundary component $\Delta_i$ satisfy $m_i=\infty$, the classical definition (c.f. \cite[Definition 9.4]{Cam11b}) asserts that $f(\BP^1) \subset X\backslash\Delta_i$. On the other hand, one can follow the definition in the arithmetic setting: the exact analogy is to fix a finite set $S$ of closed points in $\BP^1$ and to insist that the image $f(\BP^1)$ can only intersect the boundary $\Delta_i$ at these closed points. However, as shown in the next example, this definition would be again too rigid, and hence justifies the current choice.
\end{remark}

\begin{example}\label{eg after rem}
    Let $n\geqslant 3$ and $(\BP^n,\Delta)$ be a Campana absolute non-klt toric orbifold. Then $(\BP^n,\Delta)$ is not Campana rationally connected in the classical sense. Indeed, if $\vsigma$ is a set of contact orders making $(\BP^n,\Delta)$ Campana rationally connected, then the classical definition forces us to send $|\vsigma| = n+1$ many points of the curve to $\Delta$. Then the proof of \cite[Theorem 8.2]{CLT26} shows that the parameter space of the family $\pi:U\to W$ making $(\BP^n,\Delta)$ Campana rationally connected has dimension $n$. In particular, the evaluation map $ev^{(2)}$ is not dominant.
\end{example}

\section{Smooth Campana non-klt toric orbifold}\label{Smooth Campana non-klt toric orbifold}

Throughout the section, we assume that the base field $\bk$ is algebraically closed of characteristic $p>0$. 
We study the separable Campana rational connectedness of smooth non-klt Campana toric orbifolds, extending \cite[Corollary 8.4]{CLT26} (see also Corollary \ref{campana klt}), where the klt case is treated.

\begin{theorem}\label{General smooth case}
     Let $(X,\Delta_{\epsilon})$ be a Campana smooth toric orbifold such that $\mathrm{dim}(X) :=d \leqslant p$. Then it is separably Campana rationally connected. 
\end{theorem}

\begin{proof}
   We first prove the statement when $(X,\Delta_{\epsilon})$ is absolute non-klt. Choose a maximal cone $\sigma$ of the fan and reorder the primitive rays so that the first $d$ columns of the ray matrix
    \[
        R=[\rho_1\ \cdots\ \rho_n]\in M_{d\times n}(\BZ)
    \]
    are the rays of $\sigma$.
    We will prove that there exists $a_1,\dots,a_n \in \BZ_{\geqslant 0}$ such that if $L$ is the lattice generated by the $a_i\rho_i$'s, then
    \begin{equation}\label{sum}
        \sum_{i=1}^n a_i \rho_i =0,
    \end{equation}
    \begin{equation}\label{rank}
        \rank(L)=d,
    \end{equation}
    \begin{equation}\label{quotient}
        \text{and}\ \BZ^d/L \mathrm{~has~no~} p\text{-}\mathrm{torsion}.
    \end{equation}
    We will then conclude using Theorem \ref{QileMain}. Write
    \[
        R=[B\mid A],\qquad B\in\operatorname{GL}_d(\mathbb Z),\; A\in M_{d\times(n-d)}(\mathbb Z),
    \]
    so $B$ is invertible by smoothness of $\sigma$. Solutions of $Rc=0$ are parameterized by $c=(c_B,c_A) \in \BZ^{d} \times \BZ^{n-d}$ with
    \[
        B c_B + A c_A = 0\qquad\Longleftrightarrow\qquad c_B = -B^{-1}A c_A,
    \]

    Consider the induced ${\BF}_p$-linear map
    \[
        \Phi:({\BF}_p)^{n-d}\longrightarrow({\BF}_p)^d,\qquad
        \Phi(\overline{c}_A) = -\overline B^{-1}\,\overline A\,\overline{c}_A,
    \]
    and write $V=\operatorname{im}\Phi\subset({\BF}_p)^d$. It is enough to show that there exists $v \in V$ with all coordinates non-zero in $\BF_p$. Indeed, assume such $v$ exists and let $\overline{c}_A\in ({\BF}_p)^{n-d}$ such that $\Phi(\overline{c}_A)=v$. We lift $\overline{c}_A$ to any integer vector $c_A \in \BZ^{n-d}$. Set $c_B = -B^{-1}Ac_A \in \BZ^d$ and $c=(c_B,c_A) \in \BZ^n$. Then, by construction, we have $Rc=0$. In particular, if $c=(c_1,\dots,c_n)$, we have
    \[
        \sum_{i=1}^n c_i \rho_i=0.
    \]

    On the other hand, by convex geometry, since $X$ is projective, there exists $c_1',\dots,c_n'\in \BZ_{> 0}$ such that
    \[
        \sum_{i=1}^n {c_i'}\rho_i=0. 
    \]
    Set $c_i'' =c_i+ pmc_i'$ for some positive integer $m$ so that $c_i''>0$ for all $i$. Let $\vsigma$ be the set of contact order $\{c_i''\rho_i\}_{i=1}^n$. Then
    \begin{itemize}
        \item $\sum_{i=1}^n {c_i''}\rho_i=0$, so (\ref{sum}) is satisfied.
        
        \item For every $i$, we have $c_i''=c_i$ mod $p$. In particular, none of the $c_1'',\dots,c_d''$ is divisible by $p$ as it is also the case for $c_1, \dots, c_d$ by construction of $c_B=(c_1,\dots,c_d)$. 
    \end{itemize}

    Furthermore, given that $\rho_{1}, \dots, \rho_{d}$ form a $\BZ$-basis  of $N=\BZ^d$ (the ambient lattice defining the toric variety $\UX$) by hypothesis, $\det(c_{1}''\rho_{1},\dots,c_{d}''\rho_{d})=\pm\prod_{i=1}^d c_{i}''$. As none of the $c_1'',\dots,c_d''$ is divisible by $p$, 
    \[
        \rank(c_{1}''\rho_{1},\dots,c_{n}''\rho_{n})=\rank(c_{1}''\rho_{1},\dots,c_{d}''\rho_{d})=d.
    \]
    Therefore (\ref{rank}) is satisfied for $\vsigma$. Denote 
    \[
        M:=\BZ\langle c_1''\rho_1,\dots,c_d''\rho_d\rangle \subset M':=\BZ\langle c_1''\rho_1,\dots,c_n''\rho_n\rangle.
    \]
    Then
    \begin{equation}\label{indices}
        [\BZ^d:M]=[\BZ^d:M'][M':M].
    \end{equation}
    Since $\rho_1,\dots,\rho_d$ is a basis of $\BZ^d$, we have $\BZ^d/M \cong \bigoplus_{i=1}^d \BZ/c_i''\BZ$. Since $p \nmid \prod_{i=1}^d c_i''$, it follows from (\ref{indices}) that $\BZ^d/M'$ has no $p$-torsion. Therefore, (\ref{quotient}) is satisfied for $\vsigma$. We conclude that $(X,\Delta_{\epsilon})$ is separably Campana rationally connected using Lemma \ref{existence of v} and Theorem \ref{QileMain}.

    Finally, we note that the $c_i''$ may be chosen arbitrarily large. In particular, in the non-absolute non-klt case, we may add sufficiently large multiples of $p c_i'$ to each $c_i$ so that the resulting contact orders $c_i''$ along $\Delta_i$ are greater than the corresponding multiplicities $m_i$.
\end{proof}

\begin{lemma}\label{existence of v}
    With the notations of the previous theorem, $v \in V$ with all coordinates non-zero in $\BF_p$ exists.
\end{lemma}

\begin{proof}
For each $i=1,\dots,d$, let $\pi_i:\mathbb{F}_p^d \to \mathbb{F}_p$ be the projection to the $i$-th coordinate. We first show that 
\[\pi_i \circ \Phi: (\mathbb{F}_p)^{n-d} \to \mathbb{F}_p\]
is non-zero for every $i$. Set $\overline{C} = \overline{B}^{-1}\overline{A}$. The columns of $\overline{C}$ are the coordinates of the rays $\overline{\rho}_{{d+1}}, \dots, \overline{\rho}_{n}$ in the basis given by the rays $\overline{\rho}_{1}, \dots, \overline{\rho}_{d}$. Since 
\[\pi_i(\Phi(x))=-(i\mathrm{-row~of~}\overline{C})x,\]
$\pi_i \circ \Phi$ is zero if and only if the $i$-th row of $\overline{C}$ is zero.\\
Suppose for contradiction that the $i$-th row of $\overline{C}$ is zero. Then, every $\overline{\rho}_j$, with $j >d$, has zero $i$-th coordinate in the basis $\overline{\rho}_{1}, \dots, \overline{\rho}_{d}$. Equivalently, \[\overline{\rho}_{{j}}\in W_i:=\text{Span}_{\mathbb{F}_p}( \overline{\rho}_{{1}}, \dots,\hat{\overline{\rho}}_i,\dots, \overline{\rho}_{d})\]
for every $j>d$.
   Consider the codimension-one face $\eta_i:=\mathrm{cone}(\rho_{1},\dots,\hat{\rho}_{i},\dots,\rho_{d})$ of the cone $\sigma=\mathrm{cone}(\rho_{1},\dots,\rho_{d})$. Since $X$ is complete, $\eta_i$ is contained in another maximal cone
    \[
        \tau_i=\mathrm{cone}(\rho_{1},\dots,\hat{\rho}_{i},\dots,\rho_{d},\rho'_i)
    \]
    where $\rho'_i \in \{\rho_{{d+1}},\dots,\rho_{n} \}.$
    Given that $\tau_i$ is a maximal cone and $X$ is smooth,
    \[
        \{\rho_{1}, \dots, \hat{\rho}_{i}, \dots, \rho_{d},\rho'_i \}
    \]
    form a $\BZ$-basis of $N$. Therefore, after reduction mod $p$, 
     \[
        \{\overline{\rho}_{1}, \dots, \hat{\overline{\rho}}_{i}, \dots, \overline{\rho}_{d},\overline{\rho'}_i \}
    \]
    form an $\mathbb{F}_p$-basis of  $N \otimes_{\mathbb{Z}}\mathbb{F}_p$. In particular, $\overline{\rho'}_i \notin W_i$. But ${\rho'}_i$ is one of the rays $\rho_{{d+1}},\dots,\rho_{n}$ and we have already seen that our assumption implies that $\overline{\rho}_j \in W_i$ for $j>d.$ This yields a contradiction.
    
    Therefore, for each $i$, there exists $u_i \in (\mathbb{F}_p)^{n-d}$ such that $\pi_i(\Phi(u_i))\neq 0$.
  If we let $H_i=\{(x_1,\dots,x_d)\in (\BF_p)^d: x_i= 0\}$, then $\Phi(u_i) \in V$ but $\Phi(u_i) \notin H_i$. Hence $V \cap H_i$ is a proper subspace of $V$. In particular, for each $i$, the cardinality $|V \cap H_i| \leqslant p^{\mathrm{dim}_{\mathbb{F}_p}V-1}$. Hence, since all $V \cap H_i$ have a common zero and because of the assumption that $d\leqslant p$, we have
    \[|\bigcup_{i=1}^d (V\cap H_i)| \leqslant 1 + \sum_{i=1}^d (|V \cap H_i|-1) \leqslant  1+ d(p^{\mathrm{dim}_{\mathbb{F}_p}V-1} -1) < p^{\mathrm{dim}_{\mathbb{F}_p} V}=|V|.\]

    Therefore, $ V\not\subset\bigcup_{i=1}^d (V\cap H_i)$ and there exists $v=(v_1,\dots,v_d)\in V$ with $v_i\neq 0$ for all $i=1,\dots,d$. 
\end{proof}

The proof of Theorem \ref{General smooth case} shows that separable rational connectedness can be obtained under weaker assumptions than full smoothness.

\begin{corollary}\label{3 adjacent smooth}
    Let $(X,\Delta_{\epsilon})$ be a Campana toric orbifold such that $\mathrm{dim}(X) \leqslant p$. Suppose there exists a smooth maximal cone in the fan $\Sigma$ of $\UX$ such that all of its adjacent cones (that is, sharing $\mathrm{dim}(X)-1$ rays) are smooth. Then $(X,\Delta_{\epsilon})$ is separably Campana rationally connected.
\end{corollary}

\begin{remark}
    The proof of Theorem \ref{General smooth case} shows, in particular, that every projective toric variety is rationally connected (this follows from the datum of the coefficients $c_i'$ together with Theorem \ref{QileMain}). Moreover, it says that every Campana toric orbifold, whether klt or non-klt, is Campana rationally connected.
\end{remark}

We prove here the converse of \cite[Theorem 8.2]{CLT26} which gives a combinatorial sufficient condition for separable rational connectedness of toric orbifolds. Consequently, this condition is both necessary and sufficient. We then use this combinatorial characterization to construct in characteristic $p$, a Campana toric orbifold $(X,\Delta)$ $\mathrm{dim}(X) >p$ that is not separably Campana rationally connected.

\begin{prop}\label{inverse}
    Let $\UX$ be a projective toric variety of dimension $d$ with toric boundary $\Delta$. Let $\vsigma$ be a collection of positive contact orders satisfying the balancing condition. Then, $X$ is $\vsigma$-separable rationally connected if and only if:  
    \begin{itemize}
        \item The sub-lattice $\BZ\cdot \vsigma\subset N$ is of rank $d$,
        \item $N/(\BZ\cdot\vsigma)$ contains no $p$-torsion.
    \end{itemize}
\end{prop}

\begin{proof}
    We will show the \textit{only if} part since the \textit{if} part is \cite[Theorem 8.2]{CLT26}. 
      
    Let $\vsigma:=\{c_1,\dots,c_r\}$ and assume that $X$ is $\vsigma$-separably rationally connected. There is a family of non-degenerate stable log maps $\pi: \CU \to W$, $ev:  \CU \to X$ with contact orders $\vsigma$, such that $ev^{(2)}$ is dominant and separable. Let $T \subset \underline{X}$ be the dense torus and set $\CE = (ev^{(2)})^{-1}(T \times T)$. Then $ev^{(2)}_{|\CE}$ is dominant and separable. For each $w \in W$, let $\{p_1(w),\dots,p_r(w)\} \subset \mathbb{P}^1$ (with homogeneous coordinates $p_k(w):=[s_k(w):t_k(w)]$) be the set of marked points, with contact order $c_k$ at $p_k(w)$. If $m \in M=\mathrm{Hom}(N,\mathbb{Z})$, we can write
    \[
        ev_w^* \chi^m =\lambda_m(w) \prod_{k=1}^r (ts_k(w)-st_k(w))^{\langle m,c_k \rangle}
    \]
    for some $\lambda_m(w) \in \bk^{\times}$.
    
    Let $\sigma: T \times T \to T$ be the map sending $(x,y) \mapsto xy^{-1}$ and set $\Psi:= \sigma \circ ev^{(2)}_{|\CE}$. Up to replacing $\CE$ by a dense open subset, we set $q=[s:1]$ and $q'=[s':1]$. Define the morphism
    \[
        R: \CE \to \mathbb{G}_m^r, (w,q,q')\mapsto \left(\frac{s_k(w)-st_k(w)}{s_k(w)-s't_k(w)}\right)_k. 
    \]
    Since $q$ and $q'$ are sent to $T$ by $ ev_w$ while $p_k(w)$ is sent to the boundary, $R$ is well defined. Viewing $c_k \in N=\mathrm{Hom}(\mathbb{G}_m,T)$, define
    \[
        \varphi: \mathbb{G}_m^r \to T, (u_1,\dots,u_r) \mapsto \prod_{k=1}^rc_k(u_k).
    \]
    We may check (for example by applying $\chi^m$ on a basis of $M$ to both $\sigma \circ ev^{(2)}_{|\CE}$ and $\varphi \circ R$ that 
    \[
        \Psi=\sigma \circ ev^{(2)}_{|\CE}=\varphi \circ R.
    \]
    Since $ev_{|\CE}^{(2)}$ is dominant and separable, its differential is surjective at a general point and the differential of $\sigma$ is surjective everywhere because $\sigma$ is smooth. Therefore     
    \[
        d\Psi_y:T_y\CE\to T_{\Psi(y)}T 
    \]
    is surjective at a general point. This implies that $d\varphi_{R(y)}$ is surjective. Since $\varphi$ is a homomorphism of algebraic groups, its differential has constant rank. It follows that its differential at the identity is surjective. Identifying $T_{(1,\dots,1)}\mathbb{G}_m^r \cong \bk^r$ and $T_1T\cong N\otimes_{\mathbb{Z}}\bk$, the differential of $\varphi$ at the identity is given by
    \[
        d\varphi_1:\bk^r \to N\otimes_{\mathbb{Z}}\bk,(a_1,\dots,a_r)\mapsto \sum_{k=1}^r a_k(c_k\otimes 1).
    \]
    Its surjectivity means that the vectors $c_1\otimes 1,\dots,c_r\otimes 1$ span $N \otimes_{\mathbb{Z}}\bk$ as a $\bk$-vector space. Since $\bk$ has characteristic $p$, this is equivalent to $\overline{c}_1,\dots,\overline{c}_r$ span $N\otimes_{\mathbb{Z}}\mathbb{F}_p$ as an $\mathbb{F}_p$-vector space. It follows that  $\BZ\cdot \vsigma$ has rank $d$. In particular,  $N/ \BZ\cdot \vsigma$ is finite. Tensoring with $\mathbb{F}_p$ gives a right exact sequence
    \[
        \BZ\cdot \vsigma \otimes \mathbb{F}_p\to N\otimes \mathbb{F}_p \to (N/ \BZ\cdot \vsigma)\otimes \mathbb{F}_p \to 0.
    \] 
    The first map is surjective because the $\overline{c}_k$'s span $N\otimes \mathbb{F}_p$. Hence $(N/\BZ\cdot \vsigma)\otimes \mathbb{F}_p=0$. Since $N/\BZ\cdot \vsigma $ is finite, this is equivalent to $N/\BZ\cdot \vsigma$ having no $p$-torsion.
\end{proof}

\begin{example}\label{dim=4<char=2}
    We give a counter-example to Theorem \ref{General smooth case} when we do not assume that $\mathrm{dim}(X) \leqslant p$. Let $\bk$ have characteristic $2$ and define the rays in $\mathbb{Z}^4$:
    \begin{align*}
        &\rho_1=e_1,\ \rho_2=e_2,\ \rho_3=e_3,
        \rho_4= -e_1-e_2-e_3,\\
        &\rho_5= e_4,\ \rho_6=e_1+e_2-e_4.
    \end{align*}
    Let $(X,\Delta)$ be the absolute non-klt Campana toric orbifold with maximal cones induced by $(\rho_i,\rho_j,\rho_k,\rho_l)$,  where $\{i,j,k\}\subset\{1,2,3,4\}$ and $l \in \{5,6\}$. Then by \cite[Example 7.3.5]{CLS11}, we have
    \[
        X = \BP(\CO_{\BP^1}^{\oplus 2}\oplus \CO_{\BP^1}(1)^{\oplus 2}).
    \]
    Let $\vsigma=\{\gamma_1,\dots,\gamma_s\}$ be an arbitrary collection of positive contact orders of Campana type satisfying the balancing condition. Since $X$ is smooth, if $\sigma_j$ is the minimal cone containing $\gamma_j$, we have a unique expression $\gamma_j=\sum_{\rho_i \subset\sigma_j} b_{j,i}\rho_i$ with $b_{j,i}\in \mathbb{Z}_{\geqslant 0}$. Therefore, we can write $\gamma_j=\sum_{i=1}^6 b_{j,i}\rho_i$ with $b_{j,i}=0$ if $\rho_i \not\subset \sigma_j$. Set $c_i=\sum_{j=1}^s b_{j,i}$ for $i \in \{1,\dots,6\}$. Note that $\sum_{i=1}^6 c_i \rho_i=\sum_{i=1}^6\sum_{j=1}^sb_{j,i}\rho_i=\sum_{j=1}^s\gamma_j=0$. Therefore, using the notation in Theorem \ref{General smooth case}, $Rc=0$ where $c=(c_1,\dots,c_6)^T$. Hence, reducing mod $2$ gives 
    \[
        (\overline{c}_1,\dots,\overline{c}_6)=(a+b, a+b, a, a, b, b)
    \]
    for $a,b\in\BF_2$. By the absolute non-klt condition, we see that $\overline{c}_i=0 \implies \overline{b}_{j,i}=0$ for every $j$, since at most one $b_{j,i}$ is non-zero, hence every reduced contact order $\overline{\gamma}_j$ belongs to the vector space
    \[
        W = \text{Span}(\{\overline{\rho}_i\ |\ \overline{c}_i \neq 0\}).
    \]
    By varying the choice of $a$ and $b$, we see that $W$ is always a proper linear subspace of $N\otimes\mathbb{F}_2$. Therefore, by the previous proposition, $(X,\Delta)$ is not separably Campana rationally connected. 
\end{example}

\begin{remark}\label{SCRC lifts}
    Let $(X,\Delta_{\epsilon})$ be a Campana toric orbifold. Let  $\widetilde{X}\to X$ be a toric blow-up such that $\widetilde{\Sigma}$ refines $\Sigma$. If $\widetilde{\rho} \in \widetilde{\Sigma}(1)$ already appears in $\Sigma(1)$ and is equal to some $\rho$, we define $m_{\widetilde{\rho}} = m_\rho$; otherwise, we set $m_{\widetilde{\rho}}$ to be arbitrary in $\BZ_{\geqslant 1} \cup \{\infty\}$. Let
    \[
        \widetilde{\Delta}_\epsilon = \sum_{\widetilde{\rho}\in\widetilde{\Sigma}(1)} (1 - \frac{1}{m_{\widetilde{\rho}}}) D_{\widetilde{\rho}}.
    \]
    Suppose that $(X,\Delta_\epsilon)$ is separably Campana rationally connected with contact orders which are supported on the rays. Then so is $(\widetilde X,\widetilde\Delta_\epsilon)$. Indeed, we simply keep the same set of contact orders and assign zero contact along every new ray. 
\end{remark}

\begin{example}
    Every smooth projective Campana non-klt orbifold surface is separably Campana rationally connected by good contact orders. Indeed, every smooth projective toric surface is a toric blow up of $\mathbb{P}^2$ or a Hirzebruch surface. For any multiplicities on their toric boundary, these two surfaces admit good contact orders supported on their rays, making them separably Campana rationally connected. The conclusion then follows from the previous remark. On the other hand, this also follows from Theorem \ref{General smooth case}, since in this case the dimension hypothesis is automatically satisfied.
\end{example}

\section{Singular Campana non-klt toric orbifold\label{Singular Campana non-klt toric orbifold}}

\subsection{Singular non-klt Campana toric surface}

In the previous section we have seen that every smooth non-klt Campana toric surface orbifold is separably Campana rationally connected. This suggests that any failure of separable Campana rational connectedness may be related to singularities. In this section, we characterize the toric surfaces that are Campana rationally connected in terms of the configuration of their rays (cf. Proposition \ref{Criterion_Sep}). In Theorem~\ref{Crit_sing}, we further describe this property in terms of their singularities. Again, we work over an algebraically closed field $\bk$ of characteristic $p>0$.

\begin{prop}\label{Criterion_Sep}
    Let $(X,\Delta_{\epsilon})$ be a non-klt Campana toric surface orbifold.  Let $\rho_1,\dots,\rho_n$ be the primitive rays and denote by $\overline{\rho}_l$ the image of $\rho_l$ in $N\otimes_{\BZ}\BF_p$. Assume that the vectors  $\overline{\rho}_1,\dots,\overline{\rho}_n$ span $N\otimes_{\BZ}\BF_p$. Choose indices $i'$ and $j'$ such that $u := \overline{\rho}_{i'}$ and $v := \overline{\rho}_{j'}$ form an $\BF_p$-basis of $N\otimes_{\BZ}\BF_p$. Suppose that one of the following holds:
    \begin{enumerate}[\hspace{5mm}1.]
        \item There exist an index $k$ and $a,b\in\BF_p^\times$ such that $\overline{\rho}_k = au + bv$.
        
        \item There exist $s \neq i'$ and $t \neq j'$ and $a_s, a_t \in \BF_p^\times$ such that $\overline{\rho}_s = a_su$ and $\overline{\rho}_t = a_tv$.
        \end{enumerate}
    Then $(X,\Delta_{\epsilon})$ is separably Campana rationally connected.
\end{prop}

To prove the proposition, we first begin with the following lemma:

\begin{lemma}\label{reduce mod p}
    Let $(X,\Delta_{\epsilon})$ be a non-klt Campana toric surface orbifold. Let $\rho_1, \dots,\rho_n$ be the primitive rays. If there exists $(\overline{c}_1,\dots,\overline{c}_n)\in (\BF_p)^n$ such that
    \begin{enumerate}
        \item $\sum_{i=1}^n \overline{c}_i \overline{\rho}_i\equiv 0$,
        \item and indices $i'\neq j'$ such that $\overline{c}_{i'}\overline{c}_{j'}\overline{\det(\rho_{i'},\rho_{j'})} \neq 0$,
    \end{enumerate}
    then $(X,\Delta_{\epsilon})$ is separably Campana rationally connected.
\end{lemma}
 
\begin{proof} 
    We will prove that $(\overline{c}_1,\dots,\overline{c}_n)\in (\BF_p)^n$ lifts into an element $({c_1},\dots,{c_n})\in (\BZ_{\geqslant 0 })^n$ such that 
    \begin{itemize}
        \item $\sum_{i=1}^n c_i \rho_i=0$,
        \item the lattice $M$ generated by $c_1\rho_1,\dots,c_n\rho_n$ has rank $2$, and
        \item $\BZ^2/M$ has no $p$-torsion,
    \end{itemize}
    then apply Theorem \ref{QileMain}. As shown in the proof of Theorem \ref{General smooth case}, since we can choose $c_i$ arbitrary large, we can reduce to the case where $(X,\Delta_{\epsilon})$ is absolute non-klt.
    Consider the linear map:
    \begin{align*}
        R_2:   \BZ^n & \to \BZ^2\\
                v_i & \mapsto \rho_i
    \end{align*}
    (with $\{v_i\}_i$ the canonical basis of $\BZ^n$). Consider the Smith normal form of the matrix of $R_2$ where the column are given by the rays of $\UX$:
    \[
    \begin{bmatrix}
          d_1 & 0 & \dots & 0 \\
        0 & d_2 &0 &\dots 0\\
    \end{bmatrix}
    \]
    By hypothesis, $p \nmid \det(\rho_{i'},\rho_{j'})$ for some pair $(i',j')$ and by definition, $d_2$ (hence $d_1$) divides all the $2 \times 2$ minors of $R_2$, hence $\det(\rho_{i'},\rho_{j'})$. Therefore, $p \nmid d_1d_2$. Hence the reduction map
    \[
        \ker_{\BZ}(R_2) \to \ker_{\BF_p}(\overline{R}_2)
    \]
    is surjective and we can lift $(\overline{c}_1,\dots,\overline{c}_n)\in \ker_{\BF_p}(R_2)$ into an element $({c_1},\dots,{c_n})\in \ker_{\BZ}(R_2)$ with all coordinates strictly positive by Lemma \ref{surj}. In particular, we have
    \[
        \sum_{i=1}^n c_i \rho_i=0.
    \]
    Denote $M$ the lattice generated by $c_1\rho_1,\dots,c_n\rho_n$. Since there is a pair $(i',j')$ such that $c_{i'}$, $c_{j'}$, and $\det(\rho_{i'},\rho_{j'})$ are nonzero, the rank of $M$ is $2$. In addition, the assumption $\overline{c}_{i'}\overline{c}_{j'}\overline{\det(\rho_{i'},\rho_{j'})} \neq 0$ implies that $\BZ^2/M$ has no $p$-torsion. 
\end{proof}

\begin{lemma}\label{surj}
    Let the $2 \times n$ matrix 
    \[
        UR_2V=D=\begin{bmatrix}
        d_1 & 0 & \dots & 0 \\
        0 & d_2 &0 &\dots 0\\
        \end{bmatrix}
    \] 
    be the Smith normal form of the matrix of $R_2$, with $U,V \in GL(\BZ)$. If $p \nmid d_1d_2$, then the reduction map
    \[
        \ker_{\BZ}(R_2) \to\ker_{\BF_p}(\overline{R}_2) 
    \]
    is surjective. In addition, we can choose the lift with all coordinates strictly positive.
\end{lemma}

\begin{proof}
    Since $U,V \in GL(\BZ)$, we have $\det(U), \det(V)=\pm 1$. Let $\overline{y} \in \ker_{\BF_p}(\overline{R}_2)$ so that $\overline{R}_2\overline{y}= 0$. Write $\overline{y}=\overline{V}\overline{z}$. Then
    \begin{equation}\label{z}
        \overline{UR_2V}\overline{z}=\overline{UR}_2\overline{y} = 0.
    \end{equation}
    Write $\overline{z}\coloneqq(\overline{z}_1,\dots,\overline{z}_n)$. Since $d_1, d_2 \neq 0 \mathrm{~mod~} p$, we have $\overline{z}_1=\overline{z}_2= 0$ by (\ref{z}). Choose a lift $z=(z_1,\dots,z_n) \in \BZ^n$ of $\overline{z}$ such that $z_1=z_2=0$. Hence $Dz=0$. Therefore
    \[
    R_2(Vz)=U^{-1}Dz=0,
    \]
    and thus $Vz \in \ker_{\BZ}(R_2)$. Write $Vz=(a_1,\dots,a_n)$. As in the proof of Theorem \ref{General smooth case}, by completeness of $\UX$, there exists $(b_1,\dots,b_n) \in \mathbb{Z}^n_{>0}$ such that $\sum_{i=1}^n b_i\rho_i=0$. Set $c_i:=a_i +pmb_i$ for $m\gg 0$ so that all $c_i>0$. Then $y = (c_1,\dots,c_n)$ is such a lift.
\end{proof}

\begin{proof}[Proof of Proposition \ref{Criterion_Sep}] 
    It suffices to show that the hypothesis of the proposition imply the hypothesis of Lemma \ref{reduce mod p}. In case 1, we may set 
    \[
        \overline{c}_{i'} = a,\ \overline{c}_{j'} = b,\ \overline{c}_k = -1,
    \]
    and all other coefficients $\overline{c}_i = 0$. In case 2, we set
    \[
        \overline{c}_{i'} = a_s,\ \overline{c}_{j'} = a_t,\ 
        \overline{c}_s = -1,\ \overline{c}_t = -1,\ 
    \]
    and all other coefficients zero. Then $\sum_i \overline{c}_i \overline{\rho}_i = 0$ and $\overline{c}_{i'}\overline{c}_{j'}\overline{\det(\rho_{i'},\rho_{j'})} \neq 0$. Hence the hypothesis of Lemma \ref{reduce mod p} is satisfied and we are done.

\end{proof}

\subsubsection{Singularities of toric surfaces in characteristic $p$ \label{Descrp Sing}}

Our next goal is to study the property of being separably Campana rationally connected in terms of the singularities of the Campana toric surface orbifold. We first classify toric surface singularities into \textit{wild} and \textit{tame}.

Recall that all toric surface singularities are quotient singularities of the form $\frac{1}{m}(1,a)$. More specifically, if two rays $u,v$ induce a non-smooth $2$-dimensional cone which corresponds to a singular point of the form $\frac{1}{m}(1,a)$, then locally the surface is isomorphic to $\mathbb{A}^2/\mu_m$ (the geometric quotient) where $\mu_m$ acts on the diagonal weighted by $(1,a)$. 

\begin{definition}
    \label{cone singularity} We say that a quotient singularity of the form $\frac{1}{m}(1,a)$  is \textit{wild} if $p| m$. Otherwise, we say that it is \textit{tame}. We say that a non-smooth $2$-dimensional cone $\sigma$ is wildly (resp. tamely) singular if the corresponding quotient singularity is wild (resp. tame).
\end{definition} 

We can also relate the singularity of a cone to the configuration of its rays. 
If $\UX$ is a toric surface and $\rho$, $\rho'$ are two adjacent rays spanning a cone $\sigma$, then we have:
\begin{itemize}
    \item If $\det(u_{\rho},u_{\rho'})=1$, then $\sigma$ is smooth.
    \item If $p \nmid \det(u_{\rho},u_{\rho'})$ and $\det(u_{\rho},u_{\rho'})\neq 1$, then $\sigma$ is tamely singular.
    \item If $p | \det(u_{\rho},u_{\rho'})$, then $\sigma$ is wildly singular.
\end{itemize}

\begin{example}[Weighted projective surface] \label{WPS}
    Let $N = \BZ^2$ be a lattice spanned by $e_1=(1,0)$ and $e_2=(0,1)$. Let $\Sigma$ be the fan spanned by the following rays:
    \begin{align*}
        & \rho_1 = e_1 \\
        & \rho_2 = e_2 \\
        & \rho_3 = -e_1 - p e_2
    \end{align*}
    Then the toric variety corresponding to $\Sigma$ is the weighted projective plane $\BP(1,1,p)$. Note that the cone spanned by $\rho_1$ and $\rho_3$ is wildly singular, while the two other cones are smooth.

    Let $\vsigma$ be the collection of contact orders consisting of 
    \begin{align*}
        & c_1 =  e_1 \\
        & c_2 = p e_2 \\
        & c_3 = -e_1 - p e_2
    \end{align*}
    Let $N'$ be the sub-lattice $\BZ\cdot \vsigma \subset N$ and let $\Sigma'$ be be the fan induced by the the rays that the $c_i$'s span in $N'$. Since $c_3 = -c_1 - c_2$, the toric surface induced by $\Sigma'$ is $\BP^2$. The natural inclusion $N'\subset N$ induces an inclusion on the dual lattices $M\subset M'$, which induces a toric morphism $\phi: X' \to X$ given by
    \begin{align*}
        \BP^2 & \to \BP(1,1,p)\\
        (x:y:z) &\mapsto (x:y:z^p)
    \end{align*}
    This morphism is compatible with the action of $\mu_p$ on $\BP^2$ \[((x:y:z), \zeta) \mapsto (x:y:\zeta z).\] 
    Since $\bk[x,y,z]^{\mu_p}= \bk[x,y,z^p]$, one checks that this induces an isomorphism 
    \[
        \BP^2/\mu_p \cong \BP(1,1,p),
    \]
    where the left hand side is the geometric quotient. 
\end{example}

\begin{theorem}\label{Crit_sing}
    Let $(X,\Delta_{\epsilon})$ be a non-klt Campana toric surface orbifold. If every ray $\rho$ has a non-adjacent cone which is either smooth or tamely singular, then $(X,\Delta_{\epsilon})$ is separably Campana rationally connected. 
\end{theorem}

\begin{proof} 
    Let $\rho_1,\dots,\rho_n$ be the rays of $\UX$. According to Proposition \ref{Criterion_Sep}, if $(X,\Delta_{\epsilon})$ is not separably Campana rationally connected, then we are in the following two cases:
    \begin{itemize}
        \item either $\rank(\overline{\rho}_1, \dots, \overline{\rho}_n) = 1$, or
        
        \item $\rank(\overline{\rho}_1, \dots, \overline{\rho}_n) = 2$ and if $\{\overline{u},\overline{v}\} \subset \{\overline{\rho}_1, \dots, \overline{\rho}_n\}$ is an $\BF_p$-basis of $\BF_p^2$, then every other $\overline{\rho}_i$ lies in the same direction of either $\overline{u}$ or $\overline{v}$.
    \end{itemize}
    
    In the first scenario, if $\rank(\overline{\rho}_1, \dots, \overline{\rho}_n) = 1$, then for any cone spanned by two adjacent rays $\rho_i$ and $\rho_j$, we have $\det(\rho_i, \rho_j) \equiv 0  \mathrm{~mod~} p$, i.e. all $2$-dimensional cones are wildly singular.

    In the second scenario, without loss of generality, suppose $\overline{u} = \overline{\rho}_1$. Let $u$ be the ray corresponding to $\overline{u}$. Assume that every $\overline{\rho}_i$ with $i\neq 1$ lies in the direction of $\overline{v}$. Hence $\det(\overline{\rho}_i, \overline{\rho}_j) \equiv 0  \mathrm{~mod~} p$ for all $i,j \neq 1$. This suggests that if $\sigma$ is a cone of $\UX$ not containing $u$, then $\sigma$ is wildly singular. Hence in both scenarios, there exists a ray $\rho$ among the $\rho_i's$ such that any cone not containing $\rho$ is wildly singular.
\end{proof}

\begin{remark}
    If $\UX$ has more than or equal to four rays, the assumption of the previous theorem is equivalent to the existence of two non-adjacent cones in the fan of $\UX$ which are either smooth or tamely singular. If $\UX$ has strictly less than four rays, this is equivalent to no cone being wildly singular.
\end{remark}

\begin{figure}[h]\label{fig: singular surface}
    \begin{minipage}{0.45\textwidth}
        \begin{tikzpicture}[scale=0.9]
          
            \draw[->, gray] (-3,0) -- (3,0);
            \draw[->, gray] (0,-3) -- (0,3);

            \coordinate (r1)  at (1,0);
            \coordinate (r1p) at (-2,3);     
            \coordinate (r)   at (-1,1);
            \coordinate (r2)  at (-2,-1);
            \coordinate (r2p) at (2,-3);     

            \draw[thick] (0,0) -- (r1);
            \draw[thick] (0,0) -- (r1p);
            \draw[thick] (0,0) -- (r2);
            \draw[thick] (0,0) -- (r2p);
            
            \node[above right] at (r1) {$r_1=(1,0)$};
            \node[above left] at (r1p) {$r_1'=(-2,3)$};

            \node[left] at (r2) {$r_2=(-2,-1)$};
            \node[below right] at (r2p) {$r_2'=(2,-3)$};

            \fill[blue!20, opacity=0.2] (0,0) -- (r1) -- (r1p) -- cycle;
            \fill[green!20, opacity=0.2] (0,0) -- (r2) -- (r2p) -- cycle;

            \node at (-0.2,1.0) {$\sigma_1$};
            \node at (0,-1.2) {$\sigma_2$};

            \node[Violet] at (1.5,1.5)
            {$\det(r_1,r_1')=3$};
            
            \node[OliveGreen] at (-0.5,-2.5)
            {$\det(r_2,r_2')=8$};
        \end{tikzpicture}
    \end{minipage}
\hfill
    \begin{minipage}{0.45\textwidth}
        \begin{example}
            Let $(S, \Delta)$ be an absolute non-klt Campana toric surface orbifold with fan in Figure 1.
            When $\mathrm{char} \bk \notin \{2,3\}$, the cones $\sigma_1$ and $\sigma_2$ are both tamely singular and the toric orbifold $S$ is separably Campana rationally connected with contact orders given by 
            \[
                8r_1, r_1', 3r_2.
            \]
        \end{example}
    \end{minipage}
    \caption{A fan of a singular Campana toric surface orbifold $(S, \Delta)$.}
\end{figure}

\begin{remark} 
    If $\UX$ is a projective toric surface with three adjacent smooth cones, $(X, \Delta_{\epsilon})$ is separably Campana rationally connected by Theorem \ref{Crit_sing}, which we already knew from Corollary \ref{3 adjacent smooth}. 
\end{remark}

\subsection{Weighted projective space}

In this section, we prove Theorem \ref{Intro wps not sep}. We first recall the description of a weighted projective space as a toric variety. 

\begin{definition}
    We call an unordered tuple of positive integers $(q_0,\dots,q_n)$ a well-formed weight if $n\geqslant 1$ and any set of $n$ elements among the $q_i$'s are coprime. 
\end{definition}

Let $Q = (q_0,\dots,q_n)$ be a well-formed weight. Let $v_0,\dots,v_{n}$ be vectors in $\BZ^{n}$ generating the fan of $\BP^n$ (recall that $\sum_{i=0}^n v_i=0$). Set 
\[
    e_i = \frac{1}{q_i}\cdot v_i.
\]
In particular, $\sum_{i=0}^n q_i e_i=0$. Let $N$ be the lattice generated by the $e_i$'s. A weighted projective space $\BP(Q) = \BP(q_0,\dots,q_n)$ is the toric variety associated to the fan $\Sigma$ in $N_\BR$ whose cones are generated by all proper subsets of $\{e_0,\dots, e_n\}$. 

\begin{remark}
    Note that the general definition of a weighted projective space $\BP(Q)$ does not require $Q$ to be well-formed. Nevertheless, we can normalize the weights so that any two weighted projective spaces $\BP(Q)$ and $\BP(Q')$ are isomorphic if and only if the normalized weights $\overline{Q} = \overline{Q}'$, see for example \cite{Amr89}.
\end{remark}

Denote by $\rho_i$ the ray spanned by $e_i$.

\begin{theorem}\label{wps not sep}
    Let $(X,\Delta)$ be an absolute non-klt Campana toric orbifold such that $\UX$ is a weighted projective space $\mathbb{P}(Q)$, with $Q=(q_0,\dots,q_n)$ a well-formed weight, and $\Delta$ its toric boundary divisor. Then if $\carac \bk$ divides any of the $q_i$'s, $(X, \Delta)$ is not separably Campana rationally connected.
\end{theorem}

\begin{proof}
    Let $\vsigma = \{c_i\}_{i=0}^k$ be a set of contact orders making $(X,\Delta)$ separably Campana rationally connected. Since it satisfies the balancing condition and $(X,\Delta)$ is absolute non-klt, it follows that $k=n$ and $c_i \in \rho_i$ up to ordering. Suppose $c_n = \frac{b_n}{a_n} \cdot v_n$ with $\text{gcd}(a_n, b_n) = 1$, $a_n \geqslant 1$. Since $v_0+ \dots +v_n = 0$, we have
    \[
        c_n = -\frac{b_n}{a_n}\cdot (v_0 + \dots + v_{n-1}).
    \]
    The balancing condition and the fact that $v_0,\dots,v_{n-1}$ is a basis of $\BZ^n$ imply that $c_i = \frac{b_n}{a_n} \cdot v_i = \frac{q_ib_n}{a_n} \cdot e_i$ for $i = 0,\dots, n$. Since the $c_i$'s are lattice points of $N$, we must have $\frac{q_ib_n}{a_n} \in \BZ$ for all $i=0, \dots, n$, i.e. $a_n\ |\ \text{gcd}(q_0,\dots,q_n)$. Since $Q$ is well-formed, $a_n=1$. Hence we have that $c_i$ is an integral multiple of $v_i$ for each $i=0,\dots, n$. Therefore, any set of contact orders making $(X,\Delta)$ separably Campana rationally connected is of the form
    \begin{equation}\label{contact order q_i}
        \vsigma_m \coloneqq \{m\cdot v_0, \dots, m\cdot v_n\}
    \end{equation}
    for some non-zero integer $m$.
    
    It is left to justify that if $p$ divides any of the $q_j$'s, $N/\BZ\cdot\vsigma_m$ contains $p$-torsion, hence $(X,\Delta)$ is not separably Campana rationally connected by Proposition \ref{inverse}. Indeed, denote by $M$ the lattice generated by $q_0e_0, \dots, q_ne_n$ and $L = \BZ\cdot\vsigma_m$. Since $L \subset M \subset N$, we have 
    \[
        [N:L]=[N:M][M:L].
    \]
    Furthermore, since $\sum_{i=0}^n q_ie_i= \sum_{i=0}^n v_i=0$, we have
    \[
        [N:M]= \prod_{i=0}^n q_i
    \]
    by \cite[Lemma 1]{RT12}. Therefore, if $p|q_i$, then $p \vert [N:L]=|N/L|$. Since $N/L$ is a finite abelian group, it contains $p$-torsion, and we are done.
\end{proof}

\begin{example}
    Let $(\mathbb{P}(Q),\Delta)$ be an absolute non-klt Campana toric orbifold as in Theorem~\ref{wps not sep}, and let $p$ be a prime dividing one of the $q_i$'s. Then we can construct explicitly a toric blow-up 
    \[
        \widetilde{X} \longrightarrow \mathbb{P}(Q),
    \]
    with $\widetilde{\Delta}$ the toric boundary divisor of $\widetilde{X}$ and such that $(\widetilde{X},\widetilde{\Delta})$, seen as an absolute non-klt Campana toric orbifold, is separably Campana rationally connected with good contact orders. 
    There are two cases to consider:    
    \begin{itemize}
        \item[(1)] If $Q=(1,\dots,1,q_n)$, we add the ray $e_{n+1}=-e_n$ and subdivide the fan accordingly. Let $\widetilde{X}$ be the resulting toric variety and  $\widetilde{\Delta}$ its toric boundary divisor. We define the set $\vsigma$ of contact orders as follows:
        \begin{align*}
            c_i &= e_i \qquad \text{for } i\in\{0,\dots,n-1\},\\
            c_n &= (q_n+1)e_n,\\
            c_{n+1} &= -e_n = e_{n+1}.
        \end{align*}
        By Theorem~\ref{QileMain}, $(\widetilde{X},\widetilde{\Delta})$ is separably Campana rationally $\vsigma$-connected with good contact orders. For instance, when the toric variety is $\BP(1,1,p)$, this blow-up produces the Hirzebruch surface.

        \item[(2)] {Otherwise}, we add the ray $e_{n+1}=-\sum_{i=0}^n e_i$ and subdivide. In this case, we can check that the induced toric variety $\widetilde{X}$ is a non-trivial blow-up of $X$. We define the set $\vsigma$ of contact orders to be $c_i=e_i$ for all $i$, and using Theorem~\ref{QileMain}, we verify that $(\widetilde{X},\widetilde{\Delta})$ is separably Campana rationally $\vsigma$-connected with good contact orders.
    \end{itemize}
\end{example}

\bibliography{ref}

@incollection {Abr09,
    AUTHOR = {Abramovich, Dan},
     TITLE = {Birational geometry for number theorists},
 BOOKTITLE = {Arithmetic geometry},
    SERIES = {Clay Math. Proc.},
    VOLUME = {8},
     PAGES = {335--373},
 PUBLISHER = {Amer. Math. Soc., Providence, RI},
      YEAR = {2009},
      ISBN = {978-0-8218-4476-2},
   MRCLASS = {14E30 (11G35 14E15)},
  MRNUMBER = {2498065},
MRREVIEWER = {Alexandr\ V.\ Pukhlikov},
}

@incollection{katooriginal,
AUTHOR = {Kato, Kazuya},
     TITLE = {Logarithmic structures of {F}ontaine-{I}llusie},
 BOOKTITLE = {Algebraic analysis, geometry, and number theory ({B}altimore,
              {MD}, 1988)},
     PAGES = {191--224},
 PUBLISHER = {Johns Hopkins Univ. Press, Baltimore, MD},
      YEAR = {1989},
      ISBN = {0-8018-3841-X},
   MRCLASS = {14F30 (14G20)},
  MRNUMBER = {1463703},
MRREVIEWER = {Adolfo\ Quir\'os},
}

@book{ogusloggeom, 
    AUTHOR = {Ogus, Arthur},
     TITLE = {Lectures on logarithmic algebraic geometry},
    SERIES = {Cambridge Studies in Advanced Mathematics},
    VOLUME = {178},
 PUBLISHER = {Cambridge University Press, Cambridge},
      YEAR = {2018},
     PAGES = {xviii+539},
      ISBN = {978-1-107-18773-3},
   MRCLASS = {14D06 (14A20 14M25)},
  MRNUMBER = {3838359},
MRREVIEWER = {Howard\ M.\ Thompson},
       DOI = {10.1017/9781316941614},
       URL = {https://doi.org/10.1017/9781316941614},
}

@inproceedings {Amr89,
    AUTHOR = {Al Amrani, A.},
     TITLE = {Classes d'id\'eaux et groupe de {P}icard des fibr\'es
              projectifs tordus},
 BOOKTITLE = {Proceedings of {R}esearch {S}ymposium on {$K$}-{T}heory and
              its {A}pplications ({I}badan, 1987)},
   JOURNAL = {$K$-Theory},
  FJOURNAL = {$K$-Theory. An Interdisciplinary Journal for the Development,
              Application, and Influence of $K$-Theory in the Mathematical
              Sciences},
    VOLUME = {2},
      YEAR = {1989},
    NUMBER = {5},
     PAGES = {559--578},
      ISSN = {0920-3036},
   MRCLASS = {14C22},
  MRNUMBER = {999392},
MRREVIEWER = {Steven\ E.\ Landsburg},
       DOI = {10.1007/BF00535044},
       URL = {https://doi.org/10.1007/BF00535044},
}

@article {RT12,
    AUTHOR = {Rossi, M. and Terracini, L.},
     TITLE = {Linear algebra and toric data of weighted projective spaces},
   JOURNAL = {Rend. Semin. Mat. Univ. Politec. Torino},
  FJOURNAL = {Rendiconti del Seminario Matematico. Universit\`a{} e
              Politecnico Torino},
    VOLUME = {70},
      YEAR = {2012},
    NUMBER = {4},
     PAGES = {469--495},
      ISSN = {0373-1243,2704-999X},
   MRCLASS = {14M25 (52B20)},
  MRNUMBER = {3305560},
MRREVIEWER = {Lucia\ Maria\ Marino},
}

@article {Cam94,
    AUTHOR = {Campana, Fr\'ed\'eric},
     TITLE = {Remarques sur le rev\^etement universel des vari\'et\'es
              k\"ahl\'eriennes compactes},
   JOURNAL = {Bull. Soc. Math. France},
  FJOURNAL = {Bulletin de la Soci\'et\'e{} Math\'ematique de France},
    VOLUME = {122},
      YEAR = {1994},
    NUMBER = {2},
     PAGES = {255--284},
      ISSN = {0037-9484,2102-622X},
   MRCLASS = {32J27 (14E20 32E05)},
  MRNUMBER = {1273904},
MRREVIEWER = {Yoichi\ Miyaoka},
       URL = {http://www.numdam.org/item?id=BSMF_1994__122_2_255_0},
}

@article{Cam04,
     author = {Campana, Fr\'ed\'eric},
     title = {Orbifolds, special varieties and classification theory},
     journal = {Annales de l'Institut Fourier},
     pages = {499--630},
     publisher = {Association des Annales de l'Institut Fourier},
     volume = {54},
     number = {3},
     year = {2004},
     doi = {10.5802/aif.2027},
     mrnumber = {2097416},
     zbl = {1062.14014},
     language = {en},
     url = {https://www.numdam.org/articles/10.5802/aif.2027/}
}

@article {Cam11,
    AUTHOR = {Campana, Fr\'ed\'eric},
     TITLE = {Orbifoldes g\'eom\'etriques sp\'eciales et classification
              bim\'eromorphe des vari\'et\'es k\"ahl\'eriennes compactes},
   JOURNAL = {J. Inst. Math. Jussieu},
  FJOURNAL = {Journal of the Institute of Mathematics of Jussieu. JIMJ.
              Journal de l'Institut de Math\'ematiques de Jussieu},
    VOLUME = {10},
      YEAR = {2011},
    NUMBER = {4},
     PAGES = {809--934},
      ISSN = {1474-7480,1475-3030},
   MRCLASS = {32J27 (11G05 32J25 32Q45)},
  MRNUMBER = {2831280},
       DOI = {10.1017/S1474748010000101},
       URL = {https://doi.org/10.1017/S1474748010000101},
}

@incollection {Cam11b,
    AUTHOR = {Campana, Fr\'ed\'eric},
     TITLE = {Special orbifolds and birational classification: a survey},
 BOOKTITLE = {Classification of algebraic varieties},
    SERIES = {EMS Ser. Congr. Rep.},
     PAGES = {123--170},
 PUBLISHER = {Eur. Math. Soc., Z\"urich},
      YEAR = {2011},
      ISBN = {978-3-03719-007-4},
   MRCLASS = {14E05 (14-02 14E30)},
  MRNUMBER = {2779470},
MRREVIEWER = {Alexandr\ V.\ Pukhlikov},
       DOI = {10.4171/007-1/6},
       URL = {https://doi.org/10.4171/007-1/6},
}

@misc{CC21,
      title={Plane $\mathbb{A}^1$-curves on the complement of strange rational curves}, 
      author={Qile Chen and Ryan Contreras},
      year={2021},
      eprint={2108.09024},
      archivePrefix={arXiv},
      primaryClass={math.AG},
      url={https://arxiv.org/abs/2108.09024}, 
}

@article {CZ14,
    AUTHOR = {Chen, Qile and Zhu, Yi},
     TITLE = {Very free curves on {F}ano complete intersections},
   JOURNAL = {Algebr. Geom.},
  FJOURNAL = {Algebraic Geometry},
    VOLUME = {1},
      YEAR = {2014},
    NUMBER = {5},
     PAGES = {558--572},
      ISSN = {2313-1691,2214-2584},
   MRCLASS = {14M22 (14J45 14M10)},
  MRNUMBER = {3296805},
MRREVIEWER = {Dmitry\ Kerner},
       DOI = {10.14231/AG-2014-024},
       URL = {https://doi.org/10.14231/AG-2014-024},
}

@article {CZ16,
    AUTHOR = {Chen, Qile and Zhu, Yi},
     TITLE = {{$\Bbb A^1$}-connected varieties of rank one over nonclosed
              fields},
   JOURNAL = {Math. Ann.},
  FJOURNAL = {Mathematische Annalen},
    VOLUME = {364},
      YEAR = {2016},
    NUMBER = {3-4},
     PAGES = {1505--1515},
      ISSN = {0025-5831,1432-1807},
   MRCLASS = {14G05 (14M22)},
  MRNUMBER = {3466876},
MRREVIEWER = {Yong\ Hu},
       DOI = {10.1007/s00208-015-1257-1},
       URL = {https://doi.org/10.1007/s00208-015-1257-1},
}

@article {CZ19,
    AUTHOR = {Chen, Qile and Zhu, Yi},
     TITLE = {{$\Bbb A^1$}-curves on log smooth varieties},
   JOURNAL = {J. Reine Angew. Math.},
  FJOURNAL = {Journal f\"ur die Reine und Angewandte Mathematik. [Crelle's
              Journal]},
    VOLUME = {756},
      YEAR = {2019},
     PAGES = {1--35},
      ISSN = {0075-4102,1435-5345},
   MRCLASS = {14R10 (14M22)},
  MRNUMBER = {4026447},
MRREVIEWER = {Tomasz\ Pe\l ka},
       DOI = {10.1515/crelle-2017-0028},
       URL = {https://doi.org/10.1515/crelle-2017-0028},
}

@article{KatoF,
	author = {Kato, Fumiharu},
	doi = {10.1142/S0129167X0000012X},
	fjournal = {International Journal of Mathematics},
	issn = {0129-167X},
	journal = {Internat. J. Math.},
	mrclass = {14D22 (14H10)},
	mrnumber = {1754621},
	mrreviewer = {Samuel Dalalyan},
	number = {2},
	pages = {215--232},
	title = {Log smooth deformation and moduli of log smooth curves},
	url = {http://dx.doi.org/10.1142/S0129167X0000012X},
	volume = {11},
	year = {2000}}

@book {CLS11,
    AUTHOR = {Cox, David A. and Little, John B. and Schenck, Henry K.},
     TITLE = {Toric varieties},
    SERIES = {Graduate Studies in Mathematics},
    VOLUME = {124},
 PUBLISHER = {American Mathematical Society, Providence, RI},
      YEAR = {2011},
     PAGES = {xxiv+841},
      ISBN = {978-0-8218-4819-7},
   MRCLASS = {14M25 (05A15 05E45 52B12)},
  MRNUMBER = {2810322},
MRREVIEWER = {Ivan\ Arzhantsev},
       DOI = {10.1090/gsm/124},
       URL = {https://doi.org/10.1090/gsm/124},
}

@article{CLT26,
      title={Campana rational connectedness and weak approximation}, 
      author={Qile Chen and Brian Lehmann and Sho Tanimoto},
      journal = {Algebraic Geometry},
      volume  = {13},
      number  = {5},
      pages   = {705--746},
      year    = {2026},
      DOI = {10.14231/AG-2026-021}, 
}

@article {KM99,
    AUTHOR = {Keel, Se\'an and McKernan, James},
     TITLE = {Rational curves on quasi-projective surfaces},
   JOURNAL = {Mem. Amer. Math. Soc.},
  FJOURNAL = {Memoirs of the American Mathematical Society},
    VOLUME = {140},
      YEAR = {1999},
    NUMBER = {669},
     PAGES = {viii+153},
      ISSN = {0065-9266,1947-6221},
   MRCLASS = {14J26 (14B10 14E20 14E30)},
  MRNUMBER = {1610249},
MRREVIEWER = {Alessio\ Corti},
       DOI = {10.1090/memo/0669},
       URL = {https://doi.org/10.1090/memo/0669},
}

@article {KMM92c,
    AUTHOR = {Koll\'{a}r, J\'{a}nos and Miyaoka, Yoichi and Mori, Shigefumi},
     TITLE = {Rational connectedness and boundedness of {F}ano manifolds},
   JOURNAL = {J. Differential Geom.},
  FJOURNAL = {Journal of Differential Geometry},
    VOLUME = {36},
      YEAR = {1992},
    NUMBER = {3},
     PAGES = {765--779},
      ISSN = {0022-040X,1945-743X},
   MRCLASS = {14J45},
  MRNUMBER = {1189503},
MRREVIEWER = {Yuri\ G.\ Prokhorov},
       URL = {http://projecteuclid.org/euclid.jdg/1214453188},
}
\bibliographystyle{alphaurl}

\end{document}